\documentclass{article}
\usepackage{color}
\usepackage[numbers]{natbib}
\definecolor{cdarkblue}{rgb}{0.1,0.1,0.5}
\definecolor{cdarkred}{rgb}{0.5,0.1,0.1}
\definecolor{red}{rgb}{1.0,0,0}
\newtheorem{theorem}{Theorem}

\newtheorem{corollary}{Corollary}

\newtheorem{proof}{Proof}{\itshape}{\rmfamily}

\newcommand{\B}{\mathcal{B}}

\newcommand{\ben}{\begin{enumerate}}
\newcommand{\een}{\end{enumerate}}
\newcommand{\bit}{\begin{itemize}}
\newcommand{\eit}{\end{itemize}}
\newcommand{\be}{\begin{equation}}
\newcommand{\ee}{\end{equation}}
\newcommand{\bdm}{\begin{displaymath}}
\newcommand{\edm}{\end{displaymath}}
\newcommand{\bea}{\begin{eqnarray}}
\newcommand{\eea}{\end{eqnarray}}
\newcommand{\f}[1]{\fbox}

\newcommand{\realnos}{\mbox{{\bf R}}}

\newcommand{\dfrac}[2]{\displaystyle{\frac{#1}{#2}}}

\def\squareforqed{\hbox{\rlap{$\sqcap$}$\sqcup$}}
\def\qed{\ifmmode\else\unskip\quad\fi\squareforqed}

\begin{document}


\title{An invariance principle for maps with polynomial decay of
correlations}

\author{Marta Tyran-Kami\'nska,\thanks{
Institute of Mathematics, Silesian University, ul. Bankowa 14,
40-007 Katowice, POLAND, {\tt mtyran@us.edu.pl}}} \maketitle
\begin{abstract}{We give a general method of deriving statistical
limit theorems, such as the central limit theorem and its functional
version,  in the setting of ergodic measure preserving
transformations. This method is applicable in situations where the
iterates of discrete time maps display a polynomial decay of
correlations. }
\end{abstract}

\section{Introduction}

The decay of correlations in dynamical systems, or, more generally,
the rate of approach of a given initial distribution to an invariant
one, is an area of long standing interest and research. These rates
are usually described in terms of the speed at which the iterates of
a corresponding Frobenius-Perron operator, acting on a subspace of a
functional space, decay to zero. Quasi-compactness of this operator
on the space of function of bounded variation \citep{hofbauer82}
led to an exponential decay of correlations in the case of
uniformly expanding maps on the interval. Recently, a significant
body of work has been directed at an examination of
sub-exponential decay for specific families of maps 
(\cite{isola99, liveranietal, young99, bruin03}). The simplest
example is the Manneville-Pomeau map [for fixed $\gamma>0$ let
$T_\gamma:[0,1]\to[0,1]$ be given by $T_\gamma(y)=y+y^{1+\gamma}$
(mod $1$)] for which polynomial decay was demonstrated  for
H{\"o}lder continuous functions \cite{young99}.

Throughout this paper, $(Y,\B,\nu)$ denotes a probability measure
space (a measure space with $\nu(Y)=1$) and $T:Y\to Y$ a (non-invertible) measure preserving transformation.  
Thus $\nu$ is invariant for $T$  
{\it i.e.} $\nu(T^{-1}(A))=\nu(A)$ for all $A\in \B$. Recall that
$T$ is {\it ergodic} (with respect to  $\nu$) if for each $A\in\B$
with $T^{-1}(A)=A$ we have $\nu(A)\in \{0,1\}$ 
and $T$ is {\it mixing} (with respect to $\nu$) if and only if
$$
\nu(A \cap T^{-n}(B)) \to  \nu(A) \nu(B) \qquad \mbox{for
every}\quad A,B \in \B.
$$
In terms of the correlation function
$$
\mbox{Cor}(f,g\circ T^n):=\int f(y)g(T^n(y))\nu(dy)-\int
f(y)\nu(dy)\int g(y)\nu(dy)
$$
mixing is equivalent to $\mbox{Cor}(f,g\circ T^n)\to 0$ for all
$f\in L^1(Y,\B,\nu)$ and $g\in L^\infty(Y,\B,\nu)$. The transfer
operator ${\cal P}_{T,\nu}:L^1(Y,\B,\nu)\to L^1(Y,\B,\nu)$, by
definition, satisfies
$$\int {\cal P}_{T,\nu}^nf(y)g(y)\nu(dy)=\int
f(y)g(T^n(y))\nu(dy),$$ 
which leads to 
$$
|\mbox{Cor}(f,g\circ T^n)|\le ||g||_\infty ||{\cal
P}_{T,\nu}^nf-\int f(y)\nu(dy)||_1,
$$
valid for all $f\in L^1(Y,\B,\nu)$ and $g\in L^\infty(Y,\B,\nu)$,
so if one is able to estimate $||{\cal P}^nf-\int f(y)\nu(dy)||_{
L}$ for some norm $||\cdot||_{\cal L}\ge ||\cdot||_1$, then one
obtains an upper bound on $|\mbox{Cor}(f,g\circ T^n)|$ for $g\in
L^\infty$ and $f\in { L}$. This line of approach to the decay of
correlations was taken in the work cited above. A general method
of obtaining polynomial decay of the $L^1$ norm is presented in
\cite{young99}.

In this paper we address the question of the range of validity of
the central limit theorem and its functional counterpart, and
generalize results of \citet{gordin69}, \citet{keller80},
\citet{liverani}, and \citet{viana}. 
For measurable $h:Y\to\realnos$ with $\int h(y)\nu(dy)=0$, we say
that the {\it Central Limit Theorem} (CLT) holds for $h$ if the
distributions of the random variables
$\frac{1}{\sqrt{n}}\sum_{j=0}^{n-1} h\circ T^j$ converge weakly to a
normal distribution $N(0,\sigma^2)$
$$
\lim_{n\to\infty}\nu(\{y:
\sum_{j=0}^{n-1}h(T^j(y))<\sqrt{n}t\})=\frac{1}{\sqrt{2\pi}\sigma}\int_{-\infty}^t
e^{-\frac{x^2}{2\sigma^2}}dx,\quad t\in \realnos.
$$
This will be denoted by
$$
\frac{1}{\sqrt{n}}\sum_{j=0}^{n-1}h\circ T^j\to^d \sigma N(0,1).
$$ We introduce this notation because, for  $\sigma>0$, we have $N(0,\sigma^2)=\sigma N(0,1)$,
while  $\sigma N(0,1)$
is the point measure $\delta_0$ for $\sigma=0$. 
This allows us to state our results in a unified way. There will be
always a separate issue of determining whether $\sigma$ is positive
or zero.

A stronger result than the CLT is the {\it Weak Invariance
Principle}, also called a {\it Functional Central Limit Theorem
(FCLT)}.
Let $\sigma>0$ and define
the process $\{\psi_n(t),t\in[0,1]\}$ by
$$
\psi_n(t)=\frac{1}{\sigma\sqrt{n}}\sum_{j=0}^{[nt]-1}h\circ T^j
\;\;\mbox{for}\;\;t\in [0,1],\;n\ge 1
$$
(the sum from $0$ to $-1$ is set equal to $0$). 
If $\psi_n$ converges weakly to a  standard  Brownian motion $ w$ on
$[0,1]$,
then $h$ is said to satisfy the FCLT (the distributions generated
on the Skorohod space $D[0,1]$ by the $D[0,1]$-valued random
variables $\psi_n$ converge
weakly to the standard Wiener measure \cite{billingsley68}). 

One of our main results is the following
\begin{theorem}\label{t:int}
Let $T:Y\to Y$ be ergodic with respect to the invariant measure
$\nu$ and let $h\in L^2(Y,\B,\nu)$ be such that $\int
h(y)\nu(dy)=0$. If there is $\beta>\frac{1}{2}$ such that \be
\limsup_{n\to\infty}n^\beta||{\cal P}_{T,\nu}^n h||_2<\infty,
\label{c:int}\ee then the CLT and FCLT hold for $h$ provided that
$$
\sigma=\lim_{n\to\infty}\frac{||\sum_{j=0}^{n-1}h\circ
T^j||}{\sqrt{n}}>0.
$$
\end{theorem}
Many CLT results and invariance principles for maps have been
proven, {\it cf.} the survey \cite{denker89} which, in particular,
reviews the case of uniformly expanding maps on the interval; for
mixing maps the $L_1$ norm of  ${\cal P}^nh$ decay exponentially for
functions of bounded variation thus Theorem \ref{t:int} applies.
Observe that $$ ||{\cal P}_{T,\nu}^nh||_1\le ||{\cal
P}_{T,\nu}^nh||_2\le ||{\cal P}_{T,\nu}^nh||_\infty $$ for every
$h\in L^\infty(Y,\B,\nu)$.  On the other hand, if $T$ is ergodic,
then ${\cal P}_{T,\nu}$ is a contraction in every space
$L^p(Y,\B,\nu)$, $1\le p\le \infty$. Therefore \be ||{\cal
P}_{T,\nu}^nh||_2\le ||h||_\infty^{1/2}||{\cal
P}_{T,\nu}^nh||_1^{1/2}\label{e:norms} \ee for $h\in
L^\infty(Y,\B,\nu)$. Thus  Theorem \ref{t:int} is applicable when
$h\in L^\infty(Y,\B,\nu)$ and the $L_1$ norm of ${\cal P}^nh$ decays
polynomially as $n^{-\alpha}$ with $\alpha>1$. Although the CLT for
such  decay can be deduced from the result of \citet{liverani},
Theorem \ref{t:int} gives both the CLT and FCLT. To prove only the
CLT a weaker condition than Condition \ref{c:int} is sufficient (cf.
Theorem \ref{p:CLT}) while the polynomial rate is needed in the
proof of the FCLT.

Only recently  the FCLT was established by  \citet{pollicottsharp}
for H{\"o}lder continuous functions $h$ with $\int h(y)\nu(dy)=0$
and for maps $T_\gamma$ such as the Manneville-Pomeau map under
the hypothesis that $0<\gamma<\frac{1}{3}$. The CLT was proved by
\citet{young99} by establishing that the $L_1$ norm of ${\cal
P}^nh$ decays polynomially as $n^{-\alpha}$ with
$\alpha=\frac{1}{\gamma}-1$ which is greater
than $1$ exactly when $0<\gamma<\frac{1}{2}$. 
Thus our Theorem \ref{t:int} gives both the  CLT and FCLT when
$0<\gamma<\frac{1}{2}$ for the Manneville-Pomeau map.

The structure of the paper is as follows. Section \ref{s:pre}
briefly summarizes the required background and notation. In
Section \ref{s:main} we state and prove, using ideas of
\cite{maxwell00, derrienniclin01},  our main results (Theorem
\ref{p:CLT} and Theorem \ref{t:FCLT}) from which Theorem
\ref{t:int} directly follows. We also discuss the case of
$\sigma=0$. The last section contains examples of the
applicability of our abstract theorems. As our aim was to go
beyond the exponential decay of correlations, we give some
examples of transformations for which polynomial decay of
correlations has been proved.

\section{Preliminaries}\label{s:pre}
The definition of the Frobenius-Perron (transfer) operator for $T$
depends on a given $\sigma$-finite measure $\mu$ on the measure
space $(Y,\B)$ with respect to which $T$ is {nonsingular}, {\it
i.e.} $\mu(T^{-1}(A))=0$ for all $A\in\B$ with $\mu(A)=0$. This in
turn gives rise to different operators for different underlying
measures on $\B$. 
Thus if $\nu$ is invariant for $T$, then $T$ is nonsingular and the
{\it transfer operator} ${\cal P}_{T,\nu}:L^1(Y,\B,\nu)\to
L^1(Y,\B,\nu)$ is defined as follows. For any $f\in L^1(Y,\B,\nu)$,
there is a unique element ${\cal P}_{T,\nu} f$ in $L^1(Y,\B,\nu)$
such that \be \int_{A}{\cal P}_{T,\nu}
f(y)\nu(dy)=\int_{T^{-1}(A)}f(y)\nu(dy)\qquad\mbox{for }A\in\B.
\label{fpoper}\ee We are writing here ${\cal P}_{T,\nu}$ to
underline the dependence on $T$ and $\nu$.
The Koopman operator is defined by
$$
U_Tf=f\circ T
$$
for every measurable $f:Y\to\realnos$. In particular, $U_T$ is
also well defined for $f\in L^1(Y,\B,\nu)$ and is an isometry of
$L^1(Y,\B,\nu)$ into $L^1(Y,\B,\nu)$, {\it i.e.} $||U_T
f||_1=||f||_1$ for all $f\in L^1(Y,\B,\nu)$. The following
relation holds between the operators $U_T,{\cal
P}_{T,\nu}:L^1(Y,\B,\nu)\to L^1(Y,\B,\nu)$
 \be
{\cal P}_{T,\nu} U_T f=f\;\;\mbox{and}\;\;U_T {\cal P}_{T,\nu}
f=E(f|T^{-1}(\B))\label{fpcond}
 \ee for $f\in L^1(Y,\B,\nu),$ where
$E(\cdot|T^{-1}(\B)):L^1(Y,\B,\nu)\to L^1(Y,T^{-1}(\B),\nu)$
denotes the operator of conditional expectation. Since the measure
$\nu$ is finite, we have $L^p(Y,\B,\nu)\subset L^1(Y,\B,\nu)$ for
$p\ge 1$. The operator $U_T:L^p(Y,\B,\nu)\to L^p(Y,\B,\nu)$ is
also an isometry on this space. Note that if the conditional
expectation operator $E(\cdot|T^{-1}(\B)):L^1(Y,\B,\nu)\to
L^1(Y,\B,\nu)$ is restricted to $L^2(Y,\B,\nu)$, then this is the
orthogonal projection of $L^2(Y,\B,\nu)$ onto
$L^2(Y,T^{-1}(\B),\nu).$



The significance of using the transfer operator ${\cal P}_{T,\nu}$
is that it allows a unified approach to the study of statistical
properties of the transformation $T$. Extending the approach of
\citet{gordin69}, \citet{keller80}, \citet{liverani}, and
\citet{viana} we have the following

\begin{theorem} \label{CLT1} Let $(Y,\B,\nu)$ be a
probability measure space and $T:Y\to Y$ be ergodic with respect
to $\nu$. Suppose that $h\in L^2(Y,\B,\nu)$ is such that ${\cal
P}_{T,\nu}h=0$. Then the CLT and FCLT hold for $h$ provided that
$||h||_2>0$.

Moreover, for each $n\ge 1$ we have $\mbox{Cor}(h,g\circ T^n)=0$ for
all $g\in L^2(Y,\B,\nu)$ 
and
$$||{\sum_{j=0}^{n-1} h\circ T^j }||_2= {\sqrt{n}}||h||_2.
$$
%
\end{theorem}

For a direct proof of this result see \cite{mackeytyran}, where the
proof relies on the fact that the family
$$\{T^{-n+j}(\B),\frac{1}{\sqrt{n}}h\circ T^{n-j}: 1\le j\le n,
n\ge 1\}$$ is a martingale difference array for which the central
limit theorem may be proved by using the Martingale Central Limit
Theorem (cf.  \cite[Theorem 35.12]{billingsley95}) and the Birkhoff
Ergodic Theorem. If the assumption of ergodicity
appearing in Theorem \ref{CLT1} is omitted, 
then we obtain weak convergence to mixtures of normal distributions,
that is the distributions of the random variables
$\frac{1}{\sqrt{n}}\sum_{j=0}^{n-1}h\circ T^j$ converge weakly to a
distribution with a characteristic function of the form
$\varphi(r)=E(\exp(-\frac{1}{2}r^2\eta))$ where $\eta$ is such that
$\eta\circ T=\eta$ and $\int \eta(y)\nu(dy)=\int h^2(y)\nu(dy)$.
This again is a consequence of the Birkhoff Ergodic Theorem and
another version of the Martingale Central Limit Theorem due to
\citet[Corollary p. 561]{eagleson75}.

In general, for a given $h$ the equation ${\cal P}_{T,\nu}h=0$ might
not be satisfied. Then the idea is to write $h$ as a sum of two
functions in which one satisfies the assumptions of Theorem
\ref{CLT1} while the other is irrelevant for the CLT or FCLT to
hold. This is strongly connected with the property of weak
convergence which says that if two sequences differ by a sequence
converging in probability to zero and one of them is weakly
convergent then the other is weakly convergent to the same limit
\cite[Theorem 4.1]{billingsley68}. In particular, in our setting for
the CLT for $h$ to hold it is enough to show that there is a
$\tilde{h}$ satisfying the assumptions of Theorem \ref{CLT1} such
that the sequence
$(\frac{1}{\sqrt{n}}\sum_{j=0}^{n-1}(h-\tilde{h})\circ T^j)_{n\ge
1}$ is convergent in $\nu$-measure to zero. If, additionally,  the
sequence $(\frac{1}{\sqrt{n}}\max_{1\le k\le
n}|\sum_{j=0}^{k-1}(h-\tilde{h})\circ T^j|)_{n\ge 1}$ is convergent
to zero in $\nu$-measure, then the FCLT also holds for $h$.


Finally, we illustrate Theorem \ref{CLT1} with an example. The
Chebyshev maps \cite{adler64} on $[-1,1]$ are given by $$ S_N(y) =
\cos (N \arccos y), \qquad N=0,1,\cdots $$ with $S_0(y) = 1$
and $S_1(y) = y$. 

For $N \geq 2$ they are ergodic (and in fact mixing) with respect to
the measure $\nu$ with the density
$$
g_*(y) = \dfrac{1}{\pi \sqrt{1-y^2}}.
$$
For instance, for  $N = 2$ the transfer operator on
$L^1([-1,1],\B([-1,1]),\nu)$ is given by
$$ {\cal P}_{S_2,\nu}
f(y) = \frac{1}{2 }\left [ f\left(\sqrt{\frac{1}{2} y +
\frac{1}{2}}\right) + f\left(-\sqrt{\frac{1}{2} y +
\frac{1}{2}}\right) \right ].
$$

For even $N\ge 2$ and any odd function $h:[-1,1]\to\realnos$ which
is square integrable with respect to $\nu$, we have ${\cal
P}_{S_N,\nu}h=0$. We also have ${\cal P}_{S_N,\nu}h=0$ for the
function $h(y)=y$ and all $N$ (either even or odd). By Theorem
\ref{CLT1} the CLT and FCLT hold for $h$. This gives a theoretical
basis for the numerical observations of \citet{beck2001}.

\section{The main results}\label{s:main}


In this section we state and prove our main results. We start with
the following abstract theorem which gives the CLT under less
restrictive and easily verifiable assumptions when compared with the
theorem of \cite{gordin69}.  We adapt here the ideas of
\cite{maxwell00}.

\begin{theorem}\label{p:CLT}
Let $T$ be a measure-preserving transformation on the probability
space $(Y,\B,\nu)$ and let $h\in L^2(Y,\B,\nu)$ be such that $\int
h(y)\nu(dy)=0$. Suppose that \be \sum_{n=1}^\infty
n^{-\frac{3}{2}}||\sum_{k=0}^{n-1}{\cal P}_{T,\nu}^k h||_2<\infty.
\label{concltpo}\ee Then there exists $\tilde{h}\in L^2(Y,\B,\nu)$
such that ${\cal P}_{T,\nu}\tilde{h}=0$ and
$\frac{1}{\sqrt{n}}\sum_{j=0}^{n-1}(h-\tilde{h})\circ T^j$ converges
to zero in $L^2(Y,\B,\nu)$ as $n\to\infty$.

In particular, if $T$ is ergodic, then the CLT for $h$ provided that
$||\tilde{h}||_2>0$.
\end{theorem}

\begin{proof}  For $\epsilon>0$ define $
f_\epsilon=\sum_{k=1}^\infty \frac{{\cal
P}_{T,\nu}^{k-1}h}{(1+\epsilon)^k}.$ Observe that \be
f_\epsilon=\epsilon\sum_{n=1}^\infty \frac{\sum_{k=0}^{n-1}{\cal
P}_{T,\nu}^{k}h}{(1+\epsilon)^{n+1}}\label{abeld}. \ee Since
${\cal P}_{T,\nu}$ is a contraction in $L^2(Y,\B,\nu)$, we have
$f_\epsilon\in L^2(Y,\B,\nu)$ and $h=(1+\epsilon)f_\epsilon -{\cal
P}_{T,\nu}f_\epsilon$. Let us put $$h_\epsilon=
f_\epsilon-U_T{\cal P}_{T,\nu}f_\epsilon.$$ Then ${\cal
P}_{T,\nu}h_\epsilon=0$ and \be h=h_\epsilon +\epsilon f_\epsilon
+U_T{\cal P}_{T,\nu}f_\epsilon-{\cal P}_{T,\nu}f_\epsilon.
\label{decomp}\ee Now the arguments of \cite{maxwell00} apply.
Using
$$\int h_\epsilon(y)h_\delta(y)\nu(dy)=\int f_\epsilon(y)f_\delta(y)\nu(dy)-\int{\cal
P}_{T,\nu}f_\epsilon(y){\cal P}_{T,\nu}f_\delta(y)\nu(dy)$$ and
${\cal P}_{T,\nu}f_\epsilon=(1+\epsilon)f_\epsilon -h$ for any
$\epsilon,\delta>0$ we obtain \be ||h_\epsilon-h_\delta||_2^2\le
(\epsilon+\delta)(||f_\epsilon||_2^2+||f_\delta||_2^2).
\label{estheps}\ee Condition \ref{concltpo} and Equation \ref{abeld}
imply that $\sqrt{\epsilon}||f_\epsilon||_2\to 0$ as $\epsilon\to 0$
and $$\sum_{k=1}^\infty \sqrt{\delta_k}\sup_{\delta_k\le \epsilon\le
\delta_{k-1}} ||f_\epsilon||_2<\infty,$$ where $\delta_{k}=2^{-k}$
for $k\ge 0$ (\cite[Lemma 1.]{maxwell00}). Consequently,
$\tilde{h}=\lim_{\epsilon\to 0} h_\epsilon$ exists in
$L^2(Y,\B,\nu)$ and ${\cal P}_{T,\nu}\tilde{h}=0$. 
Let $\epsilon_n=2^{-j_n}$ for $n\ge 1$ where $j_n$ is the unique
integer $j$ for which $2^{j-1}\le n<2^j$. Then
$$
\sum_{k=0}^{n-1}(h-\tilde{h})\circ
T^k=\sum_{k=0}^{n-1}(h_{\epsilon_n}-\tilde{h})\circ T^k
+\epsilon_n \sum_{k=0}^{n-1}f_{\epsilon_n}\circ T^k+U_T^n{\cal
P}_{T,\nu}f_{\epsilon_n}-{\cal P}_{T,\nu}f_{\epsilon_n}
$$
by Equation \ref{decomp}. Since ${\cal
P}_{T,\nu}(h_{\epsilon_n}-\tilde{h})=0$, we have
\begin{eqnarray}
 \nonumber \frac{||\sum_{k=0}^{n-1}(h-\tilde{h})\circ T^k||_2}{\sqrt{n}} &\le &||h_{\epsilon_n}-\tilde{h}||_2 + (\epsilon_n\sqrt{n} +
\frac{2}{\sqrt{n}})||f_{\epsilon_n}||_2  \\
   &\le & ||h_{\epsilon_n}-\tilde{h}||_2 +
6\sqrt{\epsilon_n}||f_{\epsilon_n}||_2, \label{estalpha}
\end{eqnarray}
 but the right-hand
side of this inequality converges to $0$ as $n\to\infty$, which completes the proof.
\qed\end{proof}

One situation in which all of the assumptions of the preceding
theorem are met is described in the following

\begin{corollary}\label{c:eg}
Let $T$ be a measure-preserving transformation on the probability
space $(Y,\B,\nu)$ and let $h\in L^2(Y,\B,\nu)$ be such that $\int
h(y)\nu(dy)=0$. Suppose that \be \sum_{n=1}^\infty \frac{||{\cal
P}_{T,\nu}^n h||_2}{\sqrt{n}}<\infty. \label{concltpo1}\ee Then
$$
\frac{1}{\sqrt{n}}\sum_{k=0}^{n-1}h\circ T^k\to^d \sigma N(0,1)
$$
where
$$
\sigma=\lim_{n\to\infty}\frac{1}{\sqrt{n}}||\sum_{k=0}^{n-1}h\circ
T^k||_2.
$$
\end{corollary}

By imposing stronger assumptions on the growth of the norm in
Condition \ref{concltpo} we can deduce a stronger version of the
central limit theorem. Here we adapt the ideas of
\cite{derrienniclin01, derrienniclin03}. We use the standard
notation $b(n)=O(a(n))$ if $\limsup_{n\to\infty}b(n)/a(n)<\infty$.

\begin{theorem}\label{t:FCLT}Let $T$ be a measure-preserving transformation on the probability
space $(Y,\B,\nu)$ and let $h\in L^2(Y,\B,\nu)$ be such that $\int
h(y)\nu(dy)=0$. Suppose that \be ||\sum_{k=0}^{n-1}{\cal
P}_{T,\nu}^k h||_2=O(n^\alpha) \qquad\mbox{with}
\qquad\alpha<\frac{1}{2}.\label{poly}\ee Then
$\frac{1}{\sqrt{n}}\sum_{j=0}^{n-1}(h-\tilde{h})\circ T^j$ converges
to zero $\nu-$a.e and in $L^2(Y,\B,\nu)$ as $n\to\infty$.

In particular, if $T$ is ergodic, then the CLT and FCLT hold for $h$
provided that $||\tilde{h}||_2>0$.
\end{theorem}

\begin{proof} Condition \ref{poly} and Equation
\ref{abeld} imply that 
$||f_\epsilon||_2=O(\epsilon^{-\alpha})$ as $\epsilon\to 0$. 
We are going to show that \be
||\sum_{k=0}^{n-1}(h-\tilde{h})\circ
T^k||_2=O(n^\alpha).\label{dlest}\ee Since
$$
||h_{\epsilon_n}-\tilde{h}||_2\le \sum_{k=j_n+1}^\infty
||h_{\delta_k}-h_{\delta_{k-1}}||_2,
$$
we obtain the estimate
$||h_{\epsilon_n}-\tilde{h}||_2=O(n^{\alpha-1/2})$ using inequality
\ref{estheps} and the definition of $\epsilon_n$. We also have
$\sqrt{\epsilon_n}||f_{\epsilon_n}||_2=O(n^{\alpha-1/2})$ and the
desired assertion follows from Equation \ref{estalpha}. Now the
arguments of \cite{derrienniclin01} apply. By Theorem 2.17 of
\cite{derriennic01} the estimate \ref{dlest} implies that
$(h-\tilde{h})\in (I-U_T)^\beta(L^2(Y,\B,\nu))$ for
$\frac{1}{2}<\beta<1-\alpha$. Hence by Theorem 3.2(i) of
\cite{derriennic01}, with $p=\frac{1}{2}$,  we obtain
$$
\lim_{n\to\infty}\frac{1}{\sqrt{n}}\sum_{k=0}^{n-1}(h-\tilde{h})\circ
T^k=0\qquad \nu-\mbox{a.e.}
$$
and this in turn implies that
$$
\lim_{n\to\infty}\frac{1}{\sqrt{n}}\max_{0\le
n-1}|\sum_{j=0}^{k}(h-\tilde{h})\circ T^j|=0 \qquad \nu-
\mbox{a.e.},
$$
which completes the proof. \qed\end{proof}

\begin{corollary}
Let $(Y,\B,\nu)$ be a probability measure space and $T:Y\to Y$ be
ergodic with respect to $\nu$. Let $h\in L^2(Y,\B,\nu)$ be such that
$\int h(y)\nu(dy)=0$. Then 
\be||\sum_{k=0}^{n-1}{\cal P}_{T,\nu}^k h||_2=O(1)\label{c:gl}\ee if
and only if 
there exist $\tilde{h}, f\in L^2(Y,\B,\nu)$ such that ${\cal
P}_{T,\nu}\tilde{h}=0$,  $ h=\tilde{h}+f\circ T-f $.

In particular,  under Condition \ref{c:gl} the CLT and FCLT hold for
$h$ provided that $h\neq f\circ T-f$ for any $f$. 
%
\end{corollary}

\begin{proof} Since $L^2(Y,\B,\nu)$ is a reflexive Banach
space, Condition \ref{c:gl} is equivalent to $h=g-{\cal P}_{T,\nu}g$
with some $g\in L^2(Y,\B,\nu)$ (\citet[Proposition 1]{butzer71}).
First assume that $h=\tilde{h}+f\circ T-f$ with ${\cal
P}_{T,\nu}\tilde{h}=0$. By taking $g=\tilde{h}+f\circ T$ and noting
that ${\cal P}_{T,\nu}g={\cal P}_{T,\nu}U_Tf=f$ we arrive at
$g=h+{\cal P}_{T,\nu}g$ which implies  Condition \ref{c:gl}.

Now assume that Condition \ref{c:gl} holds. Let $g$ be such that
$h=g-{\cal P}_{T,\nu}g$.
Taking ${h_1}=g-U_T{\cal P}_{T,\nu}g$ and observing that ${\cal
P}_{T,\nu}{h_1}=0$, we arrive at the decomposition $$ h={h_1} +
f\circ T-f $$ where $f={\cal P}_{T,\nu}g$.  By Theorem \ref{t:FCLT}
there is $\tilde{h}$ such that ${\cal P}_{T,\nu}\tilde{h}=0$ and
$$\frac{1}{\sqrt{n}}||\sum_{j=0}^{n-1}(h-\tilde{h})\circ T^j||_2\to
0.$$  Since $\frac{1}{\sqrt{n}}||\sum_{j=0}^{n-1}(h-{h_1})\circ
T^j||_2= \frac{1}{\sqrt{n}}||h\circ T^n-h||_2\to 0$, we get
$h_1=\tilde{h}$ because ${\cal P}_{T,\nu}(h_1-\tilde{h})=0$ implies
$||\sum_{j=0}^{n-1}(h_1-\tilde{h})\circ
T^j||_2=\sqrt{n}||h_1-\tilde{h}||_2$, which completes the proof.
\qed\end{proof}

Now we give a simple result that derives a CLT and FCLT from a
decay of correlations. Although the CLT in this case is due to
\cite{liverani}, we also obtain the functional version.

\begin{corollary}\label{t:decay}
Let $(Y,\B,\nu)$ be a probability measure space, $T:Y\to Y$ be
ergodic with respect to $\nu$, and let $h\in L^\infty(Y,\B,\nu)$ be
such that $\int h(y)\nu(dy)=0$. Suppose that there are $\beta>1$ and
$c>0$ such that \be \left|\int h(y)g(T^n(y))\nu(dy)\right|\le
\frac{c}{n^\beta}||g||_\infty\label{corrdecay} \ee for all $g\in
L^\infty(Y,\B,\nu)$ and sufficiently large $n$. Then $\sigma\ge 0$
given by
$$\sigma^2=\int h^2(y)\nu(dy)+2\sum_{n=1}^\infty\int
h(y)h(T^n(y))\nu(dy)$$ is finite and if $\sigma>0$ the CLT and FCLT
hold for $h$.

Moreover, $\sigma=0$ if and only if $h=f\circ T-f$ for some $f\in
L^1(Y,\B,\nu)$.
\end{corollary}

\begin{proof} Condition \ref{corrdecay} implies that \be
||{\cal P}_{T,\nu}^n h||_2\le ||h||_{\infty}^{1/2}||{\cal
P}_{T,\nu}^n h||_1^{1/2} \qquad\mbox{and}\qquad ||{\cal P}_{T,\nu}^n
h||_1\le \frac{c}{n^{\beta}}\label{normestimate}. \ee (cf.
\cite{melbourne}, Proposition 1). Since all assumptions of Theorem
\ref{t:FCLT} are met and the series $\sum_{n=1}^\infty \int
h(y)h(T^n(y))\nu(dy)$ is convergent, the assertions follow. It
remains to discuss the case of $\sigma=0$. As in the proof of
Theorem \ref{p:CLT} let $ f_\epsilon=\sum_{k=1}^\infty \frac{{\cal
P}_{T,\nu}^{k-1}h}{(1+\epsilon)^k}$ .  Then the estimate of the norm
$||{\cal P}_{T,\nu}^n h||_1$ allows us to conclude that $f_\epsilon$
converges as $\epsilon\to 0$ to $\tilde{f}=\sum_{k=0}^{\infty}{\cal
P}_{T,\nu}^{k}h$ and $\tilde{f}\in L^1(Y,\B,\nu)$. From Equation
\ref{decomp} it then follows that $ h=U_Tf-f $ where $f={\cal
P}_{T,\nu}\tilde{f}$, which completes the proof. \qed\end{proof}

\section{Some examples}

\subsection{Maps with a neutral fixed point} Let $Y=[0,1]$ and
$\B=\B([0,1]$ be the $\sigma$-algebra of Borel subsets of $[0,1]$.
For fixed $\gamma>0$ let us consider the map $T_\gamma:[0,1]\to
[0,1]$ given by
\begin{equation}
T_\gamma(y)=\left\{
\begin{array}{ll}
y(1+2^\gamma y^\gamma) &0 \leq y \leq \frac {1}{2} \\
2y-1 &\frac {1}{2} < z \leq 1
\end{array}
\right.\label{lsvmap}
\end{equation}
which was introduced by \citet{liveranietal} to illustrate a
probabilistic approach to prove polynomial decay of correlations.
The transformation $T_\gamma$ is a simple model of maps with a
neutral (indifferent) fixed point at $p=0$, {\it i.e.}
$T_\gamma(p)=p$ and $|T_\gamma'(p)|=1$. As shown in
\cite{liveranietal} the transformation $T=T_\gamma$ has a unique
absolutely continuous invariant probability measure
$\nu=\nu_\gamma$, whose density is Lipschitz continuous on any
interval $(\epsilon, 1]$, and
 for each $h\in C^1([0,1])$
there exists a constant $C=C(h)$ such that for all $g\in
L^\infty([0,1],\B([0,1]),\nu)$ and $n\ge 1$
\begin{equation}
\left |\int h(y)g(T^n(y))\nu(dy)-\int h(y)\nu(dy)\int
g(y)\nu(dy)\right |\le C \rho_n ||g||_\infty\label{e:lsvm}
\end{equation}
where $\rho_n=n^{1-\frac{1}{\gamma}}(\log n)^{\frac{1}{\gamma}}.$

Let $0<\gamma<\frac{1}{2}$. Then there is $\beta\in
(1,\frac{1}{\gamma}-1)$ such that $\rho_n\le \frac{c_1}{n^\beta}$
for sufficiently large $n$. Thus by Corollary \ref{t:decay} the CLT
and FCLT hold for $h\in C^1([0,1])$ with $\int h(y)\nu(dy)=0$
provided that $h\neq f\circ T-f$ for any $f$.

\citet{young99} uses an abstract coupling approach to obtain
sub-exponential decay of correlations through the tail behaviour of
a return time function, 
applies her method to more general one-dimensional maps with an
indifferent fixed point, where in particular a finite number of
expanding branches are allowed and it is assumed that
$yT''(y)\approx y^\gamma$ near the indifferent fixed point, and
shows that for H{\"o}lder continuous functions $h$ on $[0,1]$ we
have $\rho_n=n^{1-\frac{1}{\gamma}}$ in Equation \ref{e:lsvm}. This
family of maps contains the interval maps with an indifferent fixed
point studied by \citet{pollicottsharp} and, in particular, the
Manneville-Pomeau map. Consequently, our Corollary \ref{t:decay}
extends Theorem 1 of \cite{pollicottsharp} to all $\gamma\in
(0,\frac{1}{2})$.

When $\gamma\in (\frac{1}{2},1)$ and  $h$ is H{\"older} continuous
with $h(0)\neq 0$ then the CLT does not hold as shown in
\cite{gouezel04b}.

\subsection{One-dimensional maps with critical points}  Consider the system
studied by \citet{bruin03}. Let $T : I \to I$ be a $C^3$ interval or
circle map with a finite set ${\cal C}$ of critical points ($c\in
{\cal C}$ if $T'(c)=0$) and no stable or neutral periodic orbit. $T$
is unimodal if it has only one critical point, and multimodal if it
has more than one. All critical points are assumed to have the same
finite critical order $l\in (1,\infty)$, {\it i.e.} for $c\in C$
there exists a diffeomorphism $ r:\realnos \to\realnos$ with
$r(0)=0$  such that for $y$ close to $c$
$$
 T(y) =\pm |r(y-c)|^l+T(c)
 $$ where the $\pm$  may
depend on the sign of $y- c$. 
For a critical point c, let $D_n(c) = |(T^n)'(T(c))|$. For
simplicity consider the case of unimodal maps. 
In \cite{bruin03} the method of \citet{young99} is adapted and the
rate $\rho_n$ in Equation \ref{e:lsvm} is related to the growth of
$D_n(c)$.  In particular, if there exists $\tilde{C}
> 0$, $\tau> 2l- 1$ such that $D_n(c) \ge \tilde{C} n^\tau$, for all $n\ge 1$, then the map $T$ has an absolutely
continuous invariant probability measure, the measure is ergodic,
and for any $\tilde{\tau} < \frac{\tau-1}{l-1}- 1$, we have $\rho_n
= n^{-\tilde{\tau}}$. Consequently, our  Corollary \ref{t:decay}
implies both the CLT and FCLT  for any H{\"o}lder continuous
function $h$.

In the study of asymptotic laws of return times in
\cite{bruinvaienti03} 
the CLT for $h=\log|T'|-\int \log|T'|(y)\nu(dy)$ is proved. It is
shown that $h\in L^2(Y,\B,\nu)$ and that the $L^2$ norm of ${\cal
P}_{T,\nu}^nh$ constitute a convergent series provided that
$D_n(c)\ge C n^\tau$ with $\tau>4l-3$ and $C>0$. Then Gordin's
theorem  as stated in \cite{viana} is used. 
Our Corollary \ref{c:eg} gives a more refined result in this case as
it can be used for  $h\in L^2(Y,\B,nu)$. Note that Theorem 1.1 of
\cite{liverani}  requires $h\in L^\infty(Y,\B,\nu)$.

\subsection{Transformations on metric spaces}
Let $Y=X$ be a  metric space with some metric $d$ and $\B=\B(X)$ be
the $\sigma$-algebra of Borel subsets of $X$. Consider a
transformation $T:X\to X$ such that $T^{-1}(x)$ is countable or
finite for each $x\in X$ and a strictly positive measurable function
$\psi:X\to\realnos$, called a potential, such that for each $x\in X$
the sum $\sum_{y\in T^{-1}(x)}\psi(y)$ is convergent. The
Ruelle-Perron-Frobenius operator is defined formally on bounded
measurable functions $\phi:X\to\realnos$ by
$$
({\cal L}_{\psi}\phi)(x)=\sum_{T(y)=x}\psi(y)\phi(y).
$$
For a thorough and up to date presentation of the concept of
Ruelle-Perron-Frobenius in studying decay of correlations we refer
to \cite{baladi}.

Recently \citet{pollicott00}, in the context of subshifts of finite
type, gave an estimate of the convergence speed of  $L^1$ norm of
iterates ${\cal L}_\psi^n \phi$, $n\ge 1$, when $\psi$ has a
summable variation. Later on, \citet{fanjiang01} extended it to
locally expansive Dini dynamical system and gave an estimate in the
supremum norm of $C(X,\realnos)$.

%


 Let us recall the setting and  notations of
\cite{fanjiang01}. Let $X$ be compact, $T$ be a continuous
transformation and $\psi$ be a continuous function. $T$ is said to
be locally expanding if there are constants $\lambda>1$ and $b>0$
such that $d(T(x),T(y))\ge \lambda$ if $d(x,y)\le b$. This implies
that $T$ is a local homeomorphism and the operator ${\cal L}_\psi$
acts on the Banach space $C(X,\realnos)$ of real valued continuous
functions equipped with the supremum norm
$||\psi||_\infty=\max_{x\in X}|\psi(x)|$. 
Recall that a
right continuous and increasing function $\omega\colon \realnos^+\to
\realnos^+$ with $\omega(0)=0$ is called a modulus of continuity.
Denote by ${\cal H}^\omega$ the space of all functions $\phi\in
C(X,\realnos)$ for which
$$[\phi]_\omega=\sup_{0<d(x,y)\le
a}{|\phi(x)-\phi(y)|\over\omega(d(x,y))} <\infty,$$ where $0<a\le b$
is a constant for which $T^{-1}(y)=\{x_1,\dots,x_n\}$ and $T$ has
local inverses $S_1,\dots,S_n$ defined on the pairwise disjoint sets
$S_j(B(y,a))$.
Finally, $\omega$ is said to satisfy the Dini
condition if $$\int_0^1{\omega(t)\over t}\,dt<\infty.$$ 


Suppose that $T$ is locally expanding and (topologically) mixing,
the modulus of continuity $\omega$ satisfies the Dini condition, and
$\psi\in{\cal H}^\omega$. From the Ruelle theorem 
proved in \cite{fanjiang01I} 
it follows that there exists  a strictly positive number $\rho$ and
a strictly positive  continuous function $\phi_*$ such that ${\cal
L}_\psi\phi_*=\rho \phi_*$, and a unique probability measure
$\mu_\psi$ such that
$$
\int {\cal L}_\psi\phi(x)\mu_\psi(dx)=\int \phi(x)\mu_\psi(dx).
$$
If we take $\phi_*$ to be normalized so $\int
\phi_*(x)\mu_\psi(dx)=1$, then for any $\phi\in C(X,\realnos)$
$$
||\rho^{-n}{\cal L}_\psi^n\phi-\phi_*\int
\psi(x)\mu_\psi(dx)||_\infty\to 0.
$$
The measure $\mu_\psi$ has the so-called Gibbs property and we call
the measure $$\nu=\phi_*\mu_\psi$$ the Gibbs measure for $T$. It is
an invariant probability measure for $T$.

Instead of working with the operator ${\cal L}_\psi$ let us consider
its normalization $\tilde{\cal L}$, which is defined as follows. Let
$$\tilde{\psi}=\psi \frac{\phi_*}{\rho \phi_*\circ T}$$ and
define
$$
\tilde{\cal L}={\cal L}_{\tilde{\psi}}.
$$
The important feature for $\tilde{\cal L}$ is that $\tilde{\cal
L}1=1$ and the transfer operator ${\cal P}_{T,\nu}$ on
$L^1(X,\B,\nu)$ and the operator $\tilde{\cal L}$ are related by
$$
{\cal P}_{T,\nu}\phi=\tilde{\cal L}\phi\qquad \nu-\mbox{a.e.}, \quad
\phi\in C(X,\realnos).
$$
This yields
$$
||{\cal P}_{T,\nu}^n\phi||_2\le ||\tilde{\cal
L}^n\phi||_\infty\quad\mbox{for}\quad  \phi\in C(X,\realnos)
$$
and Theorem 4 of \cite{fanjiang01} can be applied directly to obtain
an estimate on $||{\cal P}_{T,\nu}^n\phi||_2$ with $\int
\phi(x)\nu(dx)=0$ through the rate of decay to zero of
$||\tilde{\cal L}^n\phi||_\infty$ which depends on the modulus of
continuity of $\phi$ and the choice of  $\omega$, so that we limit
ourselves to recall two consequences of the estimates in
\cite{fanjiang01}:
\begin{enumerate}
\item
Let $\omega(t)\le C t^\theta$ for some constants $C>0$ and
$0<\theta\le 1$. Then 
${\cal H}^\omega={\cal C}^\theta$ is the space of
$\theta$-H{\"o}lder continuous functions. Thus $\psi\in {\cal
C}^\theta$ and it is known that the convergence speed is
exponential, so that there are constants $C>0$ and $\vartheta>0$
such that for any $\phi \in {\cal C}^\theta$ with $\int
\phi(x)\nu(dx)=0$
$$
||{\cal L}_{\tilde{\psi}}^n\phi||_\infty\le C e^{-\vartheta n},\quad
n\ge 1.
$$
\item\label{e:fanest}
Let $\omega(t)=\frac{1}{|\log t|^{\frac{3}{2}+\varepsilon}}$ and $\omega_0(t)=\frac{1}{|\log t|^{1+\varepsilon}}$
with $\varepsilon>0$. 
If the potential $\psi\in {\cal H}^\omega $ and $ \phi\in {\cal
H}^{\omega_0}$ with $\int \phi(x)\nu(dx)=0$, then there exists a
constant $C>0$ such that
$$
||{\cal L}_{\tilde{\psi}}^n\phi||_\infty\le C \frac{(\log
n)^{\frac{3}{2}+\varepsilon}}{n^{\frac{1}{2}+\varepsilon}},\quad
n\ge 1.
$$
\end{enumerate}
From  Theorem \ref{t:int} it follows that the CLT and FCLT hold for
$\phi$ in both cases. In the case when $\psi\in {\cal H}^{\omega_1}$
with $\omega_1(t)= \frac{1}{|\log t|^{2+\varepsilon}}$, it was
proved in \cite[Theorem 6.]{fanjiang01} that the CLT holds for
$\phi\in {\cal H}^{\omega_0}$ as in \ref{e:fanest}.  Thus Theorem
\ref{t:int} generalizes the result of \cite{fanjiang01}.

\section{Conclusions}

Here we have reviewed and extended central and functional central
limit theorems as established by particular types of temporal
decay of correlations. In particular, for the first time, we have
established criteria for CLT and FCLT validity based on polynomial
decay of correlations. Three concrete examples demonstrate the
utility of these results, and show that they are applicable
directly after establishing the decay through, for example, the
coupling method of Young \cite{young98, young99} which is very
flexible or through functional-analytic method using Ruelle's
operator. Another  method has been introduced in \cite{liverani95}
to deal with maps with discontinuities and to obtain exponential
decay. It involves a direct study of the Ruelle-Perron-Frobenius
operator but using the so-called Birkhoff metrics and the notion
of invariant cones. Moreover it has been adapted in \cite{maume01}
to deal with systems with sub-exponential decay. See the excellent
texts \cite{baladi, viana} for detailed discussions of the
functional-analytic methods.

\section*{Acknowledgments}
This work was supported by  the Natural Sciences and Engineering
Research Council (NSERC grant OGP-0036920, Canada). This  paper
was written while the author was visiting McGill University, whose
hospitality and support are gratefully acknowledged. The author
would like to thank Professor Michael C. Mackey for his interest
and valuable comments. 

\bibliographystyle{apalike}
\bibliography{zpf}
\end{document}